

%
%




\input amstex
\documentstyle{amsppt}
\NoBlackBoxes
\settabs 10 \columns

\def\qed {\ \ \ \vrule height6pt  width6pt depth0pt}
\def\L{\Cal L}
\def\N{\Bbb N}
\def\C{\Bbb C}
\def\P{\Cal P}
\def\F{\Bbb F}
\def\A{\Bbb A}

\topmatter

\title Completely bounded isomorphisms of operator algebras
\endtitle

\author Alvaro Arias
\endauthor

\thanks
Supported in part by NSF DMS 93-21369
\endthanks

\address 
Division of Mathematics, Computer Science and Statistics, 
The University of Texas at San Antonio, San
Antonio, TX 78249, U.S.A.\endaddress
\email arias\@ringer.cs.utsa.edu\endemail

\abstract 
In this paper the author proves that any two elements from one of the 
following classes of  operators are completely isomorphic to each other.

\item{1.} $\{VN(F_n):n\geq 2\}$. The $II_1$ factors generated by the left
regular representation of the free group on $n$-generators.

\item{2.} $\{C_\lambda^*(F_n):n\geq 2\}$. The reduced $C^*$-algebras of 
the free group on $n$-generators.

\item{3.} Some ``non-commutative'' analytic spaces introduced by
           G. Popescu [Po].

\noindent The paper ends with some applications to Popescu's version
of Von Neumann's inequality.
\endabstract

\subjclass Primary 47D25
Secondary 46L89
\endsubjclass

\endtopmatter

\heading
1. Introduction and preliminaries
\endheading

E. Christensen and A. M. Sinclair [CS] showed that any non-elementary
injective von Neumann algebra on a separable Hilbert space is
completely isomorphic to $B(H)$, and A. G. Robertson and S. Wassemann [RW]
generalized the work on [CS] and proved that an infinite dimensional
injective operator system on a separable Hilbert space is completely
isomorphic to either $B(H)$ or $\ell_\infty$.

The techniques on those papers depend on the injectivity of the spaces and
do not extend to interesting non-injective von Neumann algebras or
operator algebras. In the present
note we address some of these examples. 
For instance, we prove that all the von Neumann factors
of the free group on $n$ generators, $n\geq2$, are completely isomorphic
to each other. We prove the same result for the reduced $C^*$-algebras
of the free group on $n$-generators, $n\geq2$; and for some non-selfadjoint
operator algebras introduced by G. Popescu [Po].

\bigbreak
Let $H$ be a Hilbert space and $B(H)$ the set of bounded linear
operators on $H$. If we identify $M_n(B(H))$, the set of $n\times n$
matrices with entries from $B(H)$, with $B(\ell_2^n(H))$, we have a
natural norm on $M_n(B(H))$. (Here $\ell_2^n(H)$ means
$H\oplus H\oplus\cdots\oplus H$, $n$-times).

An {\sl operator space} $X$ is a closed subspace of $B(H)$. Then considering
$M_n(X)$ as a subspace of $M_n(B(H))\equiv B(\ell_2^n(H))$, we 
have norms for $M_n(X)$, $n\geq 1$. (See [BP] and [EF] for 
more on the development of this recent theory).

Let $X, Y$ be operator spaces and $u:X\to Y$ be a linear map.
Define $u_n:M_n(X)\to M_n(Y)$ by
$$u_n\bigl[(x_{ij})\bigr]=\bigl[(u(x_{ij}))\bigr].$$
We say that $u$ is completely bounded (cb in short) if
$$\|u\|_{cb}=\sup_{n\geq 1}\|u_n\|<\infty.$$
If $u_n$ is an isometry for every $n\geq 1$, then $u$ is a complete
isometry. Finally, $X$ and $Y$ are completely isomorphic if 
there exists $u:X\to Y$ such that $u$ and $u^{-1}$ are
completely bounded.

\smallbreak

Let $X\subset B(H)$, $Y\subset B(K)$ be two operator spaces.
The spatial tensor product of $X$ and $Y$, $X\otimes Y$,
is the completion 
of the algebraic tensor product of $X$ and $Y$ with the norm
induced by 
$B(H\otimes_2K)$.
With this notation, $M_n(X)=M_n\otimes X$.

\smallbreak

One of the main features of operator spaces is that the
scalars are replaced by matrices (see [E]). To see this
concretely consider $X$ a finite dimensional operator
space with basis $\{e_1,\cdots,e_n\}$. A canonical 
element in $X$ looks like $\sum_{i=1}^n a_i e_i$ for some
$a_i\in\C$; whereas a canonical element in 
$M_n(X)=M_n\otimes X$ looks like $\sum_{i=1}^n A_i\otimes e_i$,
for some $A_i\in M_n$.

\smallbreak

Two of the most important operator spaces are the row and
column Hilbert spaces. In $B(\ell_2)$ define
$C=\overline{\hbox{span}}\{e_{i1}:i\in\N\}$, the column
Hilbert space, and
$R=\overline{\hbox{span}}\{e_{1i}:i\in\N\}$, the row
Hilbert space.
Both spaces are Banach space isometric to $\ell_2$, but
have very different operator space structure. 
If $\sum_i e_{i1}\otimes T_i\in C\otimes B(H)$, then
$$\biggl\|\sum_i e_{i1}\otimes T_i\biggr\|=
\biggl\|\sum_i T_i^* T_i\biggr\|^{1\over 2},$$
and if $\sum_i e_{1i}\otimes T_i\in R\otimes B(H)$, then
$$\biggl\|\sum_i e_{1i}\otimes T_i\biggr\|=
\biggl\|\sum_i T_i T_i^*\biggr\|^{1\over 2}.$$
\medbreak

In this paper we will prove that several operator algebras
are completely isomorphic to each other. The tool that we use
is Pe\l czy\'nski's  decomposition method.

This one says that if $X$ and $Y$ are Banach
spaces such that $X$ embeds complementably into $Y$,
$Y$ embeds complementably into $X$ and
$X$ and $Y$ satisfy one of the following conditions:
\item{1.} $X\approx X\oplus X$ and $Y\approx Y\oplus Y$, or
\item{2.} $X\approx (\sum_{i=1}^\infty X)_p$, $1\leq p\leq\infty$, 

\noindent then $X$ and $Y$ are
isomorphic. Moreover, if the embeddings and projections are completely
bounded, the isomorphism is a
complete isomorphism and then $X$ and $Y$ are completely isomorphic.

We will also use a variant of condition 2.

\medbreak

The main examples in this paper will be the $C^*$-algebras generated
by the left regular representation of the free group on $n$-generators,
$\lambda:F_n\to B(\ell_2(F_n))$.

Let $F_n$ be the free group on $n$-generators, $\ell_2(F_n)$ 
the Hilbert space with orthonormal basis $\{e_x:x\in F_n\}$, and
$\L$ (or $\L_n$ to avoid confusion) the linear span of the basis; i.e.,
$$\L=\biggl\{\sum_{i=1}^k a_i e_{x_i}: k\in\N, a_i\in\C, x_i\in F_n\biggr\}.$$

$\L$ is an algebra if we multiply the elements in the natural way: i.e.,
$e_x e_y=e_{xy}$. We think of $\L$ as the Laurent polynomials on
$n$ non-commutative coordinates. 
We use two norms on $\L$: The $\|\cdot\|_2$-norm,
induced by $\ell_2(G)$, and the $\|\cdot\|$-norm defined as
$$\|p\|=\sup\{\|pq\|_2:q\in\L, \|q\|_2\leq 1\}.$$
Notice that $p\in\L$ induces a left multiplication map on $\ell_2(F_n)$
and $\|p\|$ equals the operator norm of that map. $C_\lambda^*(F_n)$
is the closure of $\L$ in the norm topology of $B(\ell_2(F_n))$,
and $VN(F_n)$ the closure of $\L$ in the strong operator topology
of $B(\ell_2(F_n))$.

The following fact is well known (see [FP], Chapter 1). 
We sketch the proof to emphasize
an argument that appears repeatedly in the paper.

\proclaim{Proposition 1} Let $n, m=2,3,\cdots,\infty$. Then 
$C^*_\lambda(F_n)$ is contained completely isometrically in $C_\lambda^*(F_m)$
and there is a completely contractive projection onto $C^*_\lambda(F_n)$.
The same is true for $VN(F_n)$ and $VN(F_m)$.
\endproclaim

\noindent{\bf Proof.} If $n\leq m$ then the formal identity 
from $C_\lambda^*(F_n)$ into $C_\lambda^*(F_m)$ is a complete 
isometry. We claim that the orthogonal projection onto $\ell_2(F_n)$ is a
complete contraction from $C_\lambda^*(F_m)$ onto $C_\lambda^*(F_n)$.
Let $p\in\L_m$ and decompose it as 
$p=r+s$ where $r\in\L_n$ and $s\in(\L_n)^\perp$. 
If $q\in\L_n$ then $rq\in\L_n$ and $sq\in (\L_n)^\perp$. Therefore,
$$\|r\|=\sup_{q\in\L_n\atop \|q\|_2\leq 1}\|rq\|_2
\leq \sup_{q\in\L_n\atop \|q\|_2\leq 1}\|pq\|_2\leq \|p\|.$$
The completely bounded part is very similar.
 
To complete the circle we need to show that $C_\lambda^*(F_\infty)$ embeds
completely complemented into $C_\lambda^*(F_2)$. Assume that
$F_2$ is generated by $a, b$ and let $G$ be the subgroup generated by
$aba^{-1}$, $a^2ba^{-2}$, $a^3ba^{-3}, \cdots$. It is easy to see that
$G$ is isomorphic to $F_\infty$, $C^*_\lambda(G)$ is completely isometric to
$C_\lambda^*(F_\infty)$ and the orthogonal projection onto $\ell_2(G)$
is a complete contraction from $C_\lambda^*(F_2)$ onto
$C^*_\lambda(G)$. 

The proof for the $VN(F_n)$'s is very similar.\qed

\medbreak

\heading
2. Isomorphisms of $C_\lambda^*(F_n)$, $n\geq 2$.
\endheading

In this section we will prove that

\proclaim{Theorem 2} $C_\lambda^*(F_n)$ is completely isomorphic to 
$C_\lambda^*(F_\infty)$ when $n=2,3,4,\cdots$.
\endproclaim

\proclaim{Theorem 3} $VN(F_n)$ is completely isomorphic to
$VN(F_\infty)$ when $n=2,3,4,\cdots$.
\endproclaim

\proclaim{Remark} It is known that the $C_\lambda^*(F_n)$'s are not
$*$-isomorphic for different $n$'s (see [PV]); however, it is still
not known if the $VN(F_n)$'s are $*$-isomorphic to each other for
$n\geq2$ (see [S], Problem 4.4.44)
\endproclaim

The proof of Theorem 2 follows from Propositions 5 and 6. 
It is simple to
go from there to Theorem 3.

We need some notation.
Divide the generators of $F_\infty$ into
$\alpha_1,\alpha_2\cdots$; $e_1,e_2,\cdots$, and denote by $F_\alpha$ the
subgroup generated by the $\alpha$'s. $F_\alpha$ is isomorphic to $F_\infty$
of course. 

Let $K=\bigcup_{j=0}^\infty e_j F_\alpha$. Denote 
$\L_K=\hbox{span}\{e_x:x\in K\}$,
and let $\ell_2(K)$ be the closure of $\L_K$ in the $\|\cdot\|_2$-norm of 
$\ell_2(F_\infty)$, 
$C^*_\lambda(K)$ the closure of $\L_K$ in the $\|\cdot\|$-norm of 
$C_\lambda^*(F_\infty)$, and
$VN(K)$ the closure of $\L_K$ in the strong operator topology of 
$B(\ell_2(F_\infty))$.

\proclaim{Proposition 4} $C^*_\lambda(K)$ is 2-cb-complemented in 
$C^*_\lambda(F_\infty)$. Moreover, the orthogonal projection
onto $\ell_2(K)$ is completely bounded from 
$C^*_\lambda(F_\infty)$ onto $C^*_\lambda(K)$.
\endproclaim

We will present the proof of Proposition 4 after the proof of Theorem 2.

\proclaim{Proposition 5} $C^*_\lambda(K)\approx C^*_\lambda(F_\infty) \approx
C^*_\lambda(F_\infty)\oplus C^*_\lambda(F_\infty)$. Moreover, the
isomorphisms are completely bounded.
\endproclaim

\noindent {\bf Proof.} Decompose $K=K_1\bigcup K_2$ where 
$K_1=\bigcup_{j=0}^\infty e_{2j} F_\alpha$, and
$K_2=\bigcup_{j=0}^\infty e_{2j+1} F_\alpha$. 
It is clear that $C^*_\lambda(K_1)$ and  $C^*_\lambda(K_2)$ 
are completely isometric to $C^*_\lambda(K)$.
Moreover, Proposition 4 applied to $K_1$ and $K_2$ implies that they are
cb-complemented in $C^*_\lambda(F_\infty)$ by the orthogonal projection.
Therefore they are complemented in $C^*_\lambda(K)$ and we have
$$C^*_\lambda(K)=C^*_\lambda(K_1)\oplus C^*_\lambda(K_2)
\approx C^*_\lambda(K)\oplus C^*_\lambda(K).$$

Similarly, decompose $K=K_3\bigcup K_4$, where $K_3=e_1F_\alpha$ and
$K_4=\bigcup_{j=2}^\infty e_jF_\alpha$, and apply the previous argument
to conclude that 
$C^*_\lambda(K)\oplus C_\lambda^*(F_\infty)\approx C^*_\lambda(K)$.

Proposition 4 tells us that 
$C_\lambda^*(F_\infty)=C^*_\lambda(K)\oplus Z$ for some $Z$. Then
$$C_\lambda^*(F_\infty)\approx C^*_\lambda(K)\oplus Z
\approx C^*_\lambda(K)\oplus C^*_\lambda(K)\oplus Z
\approx C^*_\lambda(K)\oplus C_\lambda^*(F_\infty)
\approx C^*_\lambda(K).\qed$$

\proclaim{Proposition 6} $C_\lambda^*(F_k)\approx C_\lambda^*(F_k)\oplus 
C_\lambda^*(F_k)$, for $k=2,3,\cdots$. \endproclaim

\noindent {\bf Proof.} Divide the generators of $F_\infty$ into
$\beta_1,\cdots, \beta_k; e_1,e_2,\cdots$, and denote by $F_\beta$ the
subgroup generated by the $\beta$'s; $F_\beta$ is isomorphic to $F_k$.
Let $K_\beta=\bigcup_{j=0}^\infty e_j F_\beta$.

The proof of Proposition 4 works and we get that 
$C^*_\lambda(K_\beta)$ is 2-cb complemented in 
$C_\lambda^*(F_\infty)$; hence, by Proposition 1, it is 2-cb-complemented
in $C_\lambda^*(F_k)$ also. It is clear that 
$C^*_\lambda(K_\beta)\approx C^*_\lambda(K_\beta)\oplus C^*_\lambda(K_\beta)$
and that 
$C^*_\lambda(K_\beta)\oplus C_\lambda^*(F_k)\approx C^*_\lambda(K_\beta)$. 
Hence the  proof of Proposition 5 applies and we get the result. \qed
\medbreak

We will present the proof of Theorem 2 for completeness.
This is a standard version of Pe\l czy\'nski's decomposition method.
\smallbreak

\noindent{\bf Proof of Theorem 2.} Proposition 1 tells that
$C_\lambda^*(F_k)\approx C_\lambda^*(F_\infty)\oplus Y$ for some $Y$, and
that 
$C_\lambda^*(F_\infty)\approx C_\lambda^*(F_k)\oplus Z$, for some $Z$. 
Then Propositions 5 and 6 give
$$C_\lambda^*(F_k)\approx C_\lambda^*(F_\infty)\oplus Y\approx
C_\lambda^*(F_\infty)\oplus C_\lambda^*(F_\infty)\oplus Y\approx
C_\lambda^*(F_\infty)\oplus C_\lambda^*(F_k).$$
On the other hand
$$C_\lambda^*(F_\infty)\approx
C_\lambda^*(F_k)\oplus Z\approx 
C_\lambda^*(F_k)\oplus C_\lambda^*(F_k)\oplus Z\approx
C_\lambda^*(F_k)\oplus C_\lambda^*(F_\infty).\qed$$

\medbreak

The first step for the proof of Proposition 4 is to understand how to
norm the elements in $M_n(C^*_\lambda(K))$. The typical element in $\L_K$ looks
like: $\sum_{i\leq k}e_i p_i$, where $p_i\in\L_\alpha$; i.e., 
$p_i=\sum_j a_{ij}e_{x_j}$, for some $x_j\in F_\alpha$. When we consider
operator spaces, we replace the scalars by matrices, so the canonical
element of $M_n(C^*_\lambda(K))$ looks like
$$\sum_{i\leq k} (I\otimes e_i)A_i,\quad\hbox{for some } A_i\in M_n(\L_\alpha),
\leqno{(1)}$$
$I$ is the identity in $M_n$ and 
$$A_i=\sum_j A_{ij}\otimes e_{x_j},\quad\hbox{for some } 
A_{ij}\in M_n,\hbox{ and } x_j\in F_\alpha.\leqno{(2)}$$

We use the fact (see [HP]) that, as elements of $B(\ell_2(F_\infty))$, 
$$e_i=P_ie_i+e_iP_{-i},$$
where $P_i$ is the orthogonal projection onto the set of reduced words starting 
from a positive power of $e_i$ and $P_{-i}$ is the orthogonal projection
onto the set of reduced words starting from a negative power of $e_i$.
To simplify the notation we set $(e_i)^{-1}=e_{-i}$.

We also use that $\sum_i(P_ie_i)(P_ie_i)^*=
\sum_i P_ie_ie_{-i}P_i=\sum_i P_i\leq I$, the identity on
$B(\ell_2(F_\infty))$, and if $T_i, S_i\in B(\ell_2)$ then
$\|\sum_i T_iS_i\|\leq \|\sum_i T_iT_i^*\|^{1\over 2}
\|\sum_i S_i^*S_i\|^{1\over 2}$.

We need the following technical Lemma.

\proclaim{Lemma 7} Let $x_1, x_2\in F_\alpha$, $z_1,z_2\in F_\infty$ and
suppose that (as elements of $\ell_2(F_\infty)$)
$e_{x_1}P_{-i}e_{-i}e_{z_1}=e_{x_2}P_{-j}e_{-j}e_{z_2}$.
Then either they are equal to zero, or $i=j$, $x_1=x_2$, $z_1=z_2$ and 
$e_{x_1}P_{-i}e_{-i}e_{z_1}=e_{x_1}e_{-i}e_{z_1}.$
\endproclaim

\noindent{\bf Proof.} If $z_1$ starts from $e_i$, $P_{-i}e_{-i}e_{z_1}=0$,
so we assume that the reduced words of $z_1$ and $z_2$ do not start from
$e_i$ or $e_j$ respectively. 
Then we have that $x_1e_{-i}e_{z_1}=x_2e_{-j}e_{z_2}$. Since we have 
no cancellation on the $z$'s and $e_{-i}$ and $e_{-j}$ are  the first
non-$F_\alpha$  elements of the words, then $e_i=e_j$,
$x_1=x_2$ and $z_1=z_2$.\qed

\proclaim{Proposition 8} Let $T\in M_n(\L_K)$. Then
$$
\max\biggl\{\sup_{q\in \ell_2^n(\L_\alpha)\atop \|q\|_2\leq 1}\|Tq\|_2,
\sup_{b\in\ell_2^n\atop \|b\|_2\leq 1}\|T^*b\otimes e_0\|_2\biggr\}\leq
\|T\|\leq 
\sup_{q\in \ell_2^n(\L_\alpha)\atop \|q\|_2\leq 1}\|Tq\|_2 +
\sup_{b\in\ell_2^n\atop \|b\|_2\leq 1}\|T^*b\otimes e_0\|_2.$$
\endproclaim

\noindent{\bf Proof.} The left inequality is trivially true. For the other one
take $T\in M_n(\L_K)$ as in (1)
$$T=\sum_{i\leq k} (I\otimes e_i) A_i=\sum_{i\leq k} (I\otimes P_i e_i) A_i+
\sum_{i\leq k} (I\otimes e_i P_{-i}) A_i.$$

Then
$$\eqalign{ \biggl\|\sum_{i\leq k} (I\otimes P_i e_i) A_i\biggr\|
&\leq\biggl\|\sum_{i\leq k} (I\otimes P_i e_i)
        (I\otimes P_ie_i)^*\biggr\|^{1\over 2}
                        \biggl\|\sum_{i\leq k} A_i^* A_i\biggr\|^{1\over 2}\cr
&\leq \biggl\|\sum_{i\leq k}A_i^*A_i\biggr\|^{1\over 2}\cr
&=\sup_{q\in\ell_2^n(\L_\alpha)\atop \|q\|_2\leq 1}
     \sqrt{\sum_{i\leq k}\|A_iq\|_2^2}\cr
&=\sup_{q\in\ell_2^n(\L_\alpha)\atop \|q\|_2\leq 1}
     \biggl\|\sum_{i\leq k}(I\otimes e_i)A_iq\biggr\|_2\cr
&=\sup_{q\in\ell_2^n(\L_\alpha)\atop \|q\|_2\leq 1} \|Tq\|_2.\cr}$$

\medbreak
On the other hand, $\|\sum_{i\leq k} (I\otimes e_i P_{-i}) A_i\|=
\|\sum_{i\leq k} A_i^*(I\otimes P_{-i}e_{-i})\|$. To norm the latter
one, take $q\in\ell_2^n(\ell_2(F_\infty))$, $\|q\|_2\leq 1$ and decompose it
as
$$q=\sum_lb_l\otimes e_{z_l}\quad\hbox{ where } 
b_l\in\ell_2^n, z_l\in F_\infty.\leqno{(3)}$$

Using (2) and (3) we have
$$\sum_{i\leq k}A_i^*(I\otimes P_{-i}e_{-i}) q =
\sum_{i\leq k}\sum_j\sum_l A_{ij}^*b_l\otimes e_{-x_j}P_{-i}e_{-i}e_{z_l}.$$
Lemma 7 tells us that all those terms are orthogonal to each other 
or zero. Hence,
$$\eqalign{\biggl\|\sum_{i\leq k}A_i^*(I\otimes P_{-i}e_{-i})q\biggr\|_2^2
&\leq \sum_{i\leq k}\sum_j\sum_l\|A_{ij}^*b_l\|_2^2\cr
&= \sum_l\left[\sum_{i\leq k}\sum_j
     \biggl\|A_{ij}^*{b_l\over \|b_l\|_2}\biggr\|_2^2\right]\|b_l\|_2^2\cr
&\leq \sup_{b\in\ell_2^n\atop\|b\|_2\leq 1}
       \sum_{i\leq k}\sum_j\|A_{ij}^*b\|_2^2\cr
&=\sup_{b\in\ell_2^n\atop \|b\|_2\leq 1}\|T^*b\otimes e_0\|_2^2.\qed\cr}$$

\medbreak

\noindent{\bf Proof of Proposition 4}. Let $T\in M_n(\L_\infty)$. Write it
as $T=T_1+T_2$, where $T_1\in M_n(\L_K)$ and $T_2\in M_n((\L_K)^\perp)$. 
Notice that if $q\in\ell_2^n(\L_\alpha)$, then 
$$T_1q\in\ell_2^n(\L_K),\quad\hbox{ and }\quad T_2 q\in\ell_2^n((\L_K)^\perp).$$
Hence,
$$\sup_{q\in \ell_2^n(\L_\alpha)\atop \|q\|_2\leq 1}\|T_1q\|_2
\leq \sup_{q\in \ell_2^n(\L_\alpha)\atop \|q\|_2\leq 1}\|Tq\|_2\leq\|T\|.$$
Moreover, it is clear that given $b\in\ell_2^n$, we have that
$$\|T_1^*b\otimes e_0\|_2\leq \|T^*b\otimes e_0\|_2\leq\|T^*\|=\|T\|.$$
Therefore, by Proposition 8, $\|T_1\|\leq 2\|T\|.$\qed
\medbreak

Propositions 5 and 8 give a representation of $C_\lambda^*(F_\infty)$
in terms of row and column Hilbert spaces.

Let $T=\sum_{i\leq k}(I\otimes e_i)A_i\in M_n(\L_K)$, 
$A_i\in M_n(\L_\alpha)$. Use (2) to write
$T=\sum_{i\leq k}\sum_j A_{ij}\otimes e_ie_{x_j}$, for some
$A_{ij}\in M_n$. Then we have
$$\eqalign{
\sup_{q\in \ell_2^n(\L_\alpha)\atop \|q\|_2\leq 1}\|Tq\|_2&=
    \biggl\|\sum_{i\leq k} A_i^*A_i\biggr\|^{1\over2},\cr
\sup_{b\in\ell_2^n\atop \|b\|_2\leq 1}\|T^*b\otimes e_0\|_2&=
    \biggl\|\sum_{i\leq k}\sum_j A_{ij}A_{ij}^*\biggr\|^{1\over2}.\cr}$$
We see that the first term is the norm of $T$ in 
$M_n(\,C\otimes C_\lambda^*(\L_\alpha)\,)$, and the second one is
the norm of $T$ in $M_n(R\otimes R(F_\alpha)\,)$. Here
$R(F_\alpha)$ is $\ell_2(F_\alpha)$ with the row operator space structure.

Using the notation of interpolation theory of operator spaces (see [P]) 
we conclude

\proclaim{Proposition 9} $C_\lambda^*(F_\infty)
~{\buildrel c.b. \over \approx}~
[\,C\otimes C_\lambda^*(F_\infty)\,]~\bigcap~
[\,R\otimes R(F_\infty)\,].$
\endproclaim

\proclaim{Remark} If we are interested only in the Banach space structure, we
have that $C_\lambda^*(F_\infty)$ is isomorphic (but probably not completely
isomorphic) to $C\otimes C_\lambda^*(F_\infty)$.
\endproclaim

\heading
3. Isomorphisms of non-commutative analytic algebras. 
\endheading

In this section we will consider only the words consisting of positive
powers of the generators of $F_k$, and the identity. We denote this
set by $P_k\subset F_k$, and let $\ell_2(P_k)$ be the Hilbert space with
orthonormal basis $\{e_x:x\in P_k\}$. This Hilbert space is also
denoted by $\F^2(H_k)$, the full Fock space on $k$-generators, see [Po].
Let $\P$ (or $\P_k$ to avoid confusion) be the linear span of the basic
elements, and consider two norms on $\P$: 
The $\|\cdot\|_2$-norm,
induced by $\ell_2(P_k)$, and the $\|\cdot\|_\infty$-norm,
defined as
$$\|p\|_\infty=\sup\{\|pq\|_2:q\in\P_k, \|q\|_2\leq 1\}.$$
Notice that
$p\in\P_k$ induces a left multiplication operator on
$\ell_2(P_k)$ and $\|p\|_\infty$ equals the operator norm of
that map. 

\proclaim{Remark} If $p\in\P_k$, then the $\|p\|_\infty$-norm
does not coincide with the $\|p\|$-norm as an element of $C_\lambda^*(F_k)$.
In fact, if $Q$ is the orthogonal projection onto $\ell_2(P_k)$ then
$\|p\|_\infty=\|QpQ\|$. We always have that $\|p\|_\infty\leq \|p\|$ and
sometimes the inequality is strict (see Proposition 17).
\endproclaim

Let $\A(k)$ be the closure of $\P_k$ in the norm topology
of $B(\ell_2(P_k))$, and $\F^\infty(k)$ the closure in the 
strong operator topology. These spaces are studied in [Po], where
he calls them non-commutative analogues of the disk algebra and $H^\infty$.
When $k=1$ they coincide with the classical definitions.

The main results of this section are.

\proclaim{Theorem 10} $\A(k)$ is completely isomorphic to
$\A(\infty)$ when $n=2,3,4,\cdots.$
\endproclaim

\proclaim{Theorem 11} $\F^\infty(k)$ is completely isomorphic to 
$\F^\infty(\infty)$ when $n=2,3,4,\cdots.$
\endproclaim

As in the previous section we will only present the
proof of the first one, the other one is essentially the same.
The proof of Theorem 10 will follows from Propositions 12, 13 and 14.

\proclaim{Proposition 12} Let $n\leq m$, and let
$\Phi:\A(n)\to\A(m)$ be the formal identity. Then $\Phi$
is a complete isometry. Moreover, the orthogonal projection
onto $\ell_2(P_n)$ is completely contractive from 
$\A(m)$ onto $\A(n)$.
\endproclaim

\noindent {\bf Proof.} Let $E\subset P_m$ be the set of all $y\in P_m$
whose first letter does not start from $e_1,\cdots, e_n$. It is easy to
see that $\{P_n y:y\in E\}$ forms a partition of $P_m$. Then
$\ell_2(P_m)=\sum_{j=1}^\infty\oplus\ell_2(P_n y_j)$, where
$E=\{y_j:j\in{\N}\}$.

Let $p\in\P_n$ and $q\in\P_m$, $\|q\|_2\leq 1$. Use the previous partition
to decompose $q$ as $q=\sum_j r_je_{y_j}$, for some $r_j\in\P_n$ and 
$y_j\in E$. Then
$$\|pq\|_2^2=\sum_{i=0}^\infty \|pr_i\|_2^2
=\sum_{i=0}^\infty\biggl\|p{r_i\over \|r_i\|_2}\biggr\|_2^2\|r_i\|_2^2
\leq \sup_i \left\|p {r_i\over\|r_i\|_2}\right\|_2^2\leq\|p\|_\infty^2.$$
This tells us that $\|p\|_{\A(n)}= \|p\|_{\A(m)}$. Moreover, if 
$p\in\A(n)\subset \A(m)$, we norm it with elements from $\P_n$.

This is the fact that we use for the complementation. Given $p\in\P_m$, write
it as $p=p_1+p_2$, where $p_1\in\P_n$ and $p_2\in(\P_n)^\perp$. Then if
$q\in\P_n$, we have that $p_1q\in\P_n$ and $p_2q\in(\P_n)^\perp$. Hence,
$$\|p_1\|_\infty=\sup_{q\in\P_n\atop \|q\|_2}\|p_1q\|_2
\leq \sup_{q\in\P_n\atop \|q\|_2}\|pq\|_2\leq\|p\|_\infty.$$

The completely bounded part is very similar.\qed

\proclaim{Proposition 13} There exists a subspace of $\A(2)$
completely isometric to $\A(\infty)$ and completely
complemented by the orthogonal projection.
\endproclaim

\noindent{\bf Proof.} Let $a,b$ be the generators of $P_2$. Let $P_\alpha$
be the set of all words generated by $ab, a^2b, a^3b, a^4b,\cdots$. $P_\alpha$
is clearly isomorphic to $P_\infty$. Let $E\subset P_2$ be the set of all
words in $P_2$ such that no initial segment belongs to $P_\alpha$. Then
it is easy to see that $\{P_\alpha y:y\in E\}$ forms a partition of 
$P_2$, and the proof is like that of the previous proposition.\qed

\proclaim{Proposition 14} $\A(\infty)$ is completely isomorphic to 
$C\otimes \A(\infty)$, where $C$ is the column Hilbert space.
\endproclaim

\noindent{\bf Proof.} Divide the generators of $P_\infty$ into
$\alpha_1, \alpha_2\cdots$; $e_1, e_2,\cdots$, and let $P_\alpha$
be the set of words generated by the $\alpha$'s.  Let
$K=\bigcup_{j=1}^\infty e_j P_\alpha\subset P_\infty$, and let
$\P(K)$ be the span of the basic elements in $K$. Denote the closure
of $\P(K)$ in the $\ell_2$-norm by $\ell_2(K)$, and in the 
$\|\cdot\|_\infty$-norm 
by $\A(K)$. 

The canonical element of $\P(K)$ looks like $\sum_{i\leq k} e_ip_i$,
for some $k\in\N$ and $p_\alpha\in\P_\alpha$.
Given any $q\in\P$ the $e_ip_iq$'s are orthogonal. Then we have,
$$\biggl\|\sum_{i\leq k} e_ip_i\biggr\|_\infty
=\sup_{\|q\|_2\leq 1}\biggl\|\sum_{i\leq k} e_ip_iq\biggr\|_2
=\sup_{\|q\|_2\leq 1}\sqrt{\sum_{i\leq k}\|p_iq\|_2^2}
=\biggl\|\sum_{i\leq k} p_i^*p_i\biggr\|_\infty^{1\over2}.$$
Since $\A(P_\alpha)$ is isometric to $\A(\infty)$ 
we conclude that $\A(K)$ is isometric to $C\otimes \A(\infty)$. 
Moreover, the elements in $\A(k)$ are normed by elements in $\ell_2(P_\alpha)$.

We will now see that $\A(K)$ is complemented in $\A(\infty)$.
Let $p\in\P_\infty$ and decompose it as $p=p_1+p_2$, where
$p_1\in\P(K)$ and $p_2\in(\P(K))^\perp$. If $q\in\P_\alpha$, then 
$p_1q\in\P(K)$ and $p_2q\in(\P(K))^\perp$. Hence,
$\|p_1q\|_2\leq\|pq\|_2$, and $\|p_1\|_\infty\leq\|p\|_\infty$.

As in the proof of Proposition 5 
it is clear that $\A(K)$ is isomorphic to its square, and then 
isomorphic to $\A(\infty)$. 

The completely bounded
part follows in the same way after we replace the scalars by matrices.\qed
\medbreak

\noindent{\bf Proof of Theorem 10.}  By Proposition 12 and 13 we have that
$\A(k)=\A(\infty)\oplus Y$ for some $Y$. Since
$\A(\infty)\approx \A(\infty)\oplus \A(\infty)$, we have
$$\A(k)=\A(\infty)\oplus Y\approx \A(\infty)\oplus \A(\infty)
\oplus Y\approx \A(\infty)\oplus \A(k).$$
On the other hand, $\A(\infty)=\A(k)\oplus Z$, for some $Z$. If
$Q:\A(\infty)\to\A(k)$ is that projection and $I:C\to C$ is the
identity on $C$, $I\otimes Q$ decomposes
$C\otimes \A(\infty)=[C\otimes \A(k)]\oplus [C\otimes Z]$
because $Q$ is completely bounded. Hence,
$$\eqalign{ \A(\infty)&\approx C\otimes \A(\infty)\cr
   &=[C\otimes \A(k)]\oplus [C\otimes Z]\cr
   &\approx \A(k)\oplus [C\otimes \A(k)]\oplus [C\otimes Z]\cr
   &\approx \A(k)\oplus [C\otimes \A(\infty)]\approx 
              \A(k)\oplus\A(\infty).\qed \cr}$$

\heading
4. Applications to Von Neumann's inequality
\endheading

Fix $k$ for this section and let $P_k$ be the positive words generated
by $e_1,\cdots,e_k$. As in the previous section,
$\F^2(H_k)=\ell_2(P_k)$ is the full Fock space on $H_k$, a $k$-dimensional
Hilbert space; $\A(k)$ and $\F^\infty(k)$ have the same meaning.

In [Po] G. Popescu proved that if $T_1, \cdots, T_k\in B(\ell_2)$ are such
that $\|\sum_{i\leq k} T_iT_i^*\|\leq 1$ (i.e., $\|[T_1\cdots T_k]\|\leq1$)
then any $p(e_1,\cdots, e_k)\in\A(k)$ satisfies
$$\|p(T_1,\cdots, T_k)\|\leq \|p\|_{\A(k)}.\leqno{(4)}$$
When $k=1$, this is the classical Von Neumann's inequality.

In this section we prove that $\A(k)$ and $\F^\infty(k)$ contain many
complemented Hilbertian subspaces. Hence we can easily compute
$\|p\|_{\A(k)}$ whenever $p$ belongs to one of those subspaces,
and use $\|p\|_{\A(k)}$ in Popescu's inequality (4). (See [AP] for
more examples and connections with inner functions).

We start with the following elementary lemma.

\proclaim{Lemma 15} Let $p=\sum_i a_i e_{x_i}, q=\sum_j b_je_{y_j}\in\P$ be
such that $x_iy_j=x_{i'}y_{j'}$ if and only if $x_i=x_{i'}$ and
$y_j=y_{j'}$, (that is, we cannot have any cancellation), then
$\|pq\|_2=\|p\|_2\|q\|_2$.
\endproclaim

\noindent{\bf Proof.} We have that $pq=\sum_i\sum_ja_ib_je_{x_iy_j}$. Since
all the $x_iy_j$-terms are different, the $e_{x_iy_j}$'s are orthogonal.
Hence,
$$\|pq\|_2=\sqrt{\sum_i\sum_j|a_ib_j|^2}=\sqrt{\sum_i|a_i|^2}
\sqrt{\sum_j|b_j|^2}=\|p\|_2\|q\|_2.\qed$$

\proclaim{Remark} The lemma extends to the $M_n$-case just as easily.
If $T=\sum_i A_i\otimes e_{x_i}\in M_n(\P)$ and 
$q=\sum_j b_j\otimes e_{y_j}\in\ell_2^n(\P)$, then
$Tq=\sum_i\sum_jA_ib_j\otimes e_{x_iy_j}$ and
$\|Tq\|_2=\sqrt{\sum_i\sum_j \|A_ib_j\|_2^2}$.
\endproclaim

\proclaim{Proposition 16} Let $W_n\subset P_k$ be the set of all words
in $P_k$ having $n$-letters, and let 
$X_n=\hbox{span}\{e_x:x\in W_n\}\subset \A(k)$.
Then $X_n$ is completely isometric to $C_{n^k}$, the column Hilbert space
of the same dimension. Moreover, $X_n$ is completely complemented in $\A(k)$.
\endproclaim

\noindent{\bf Proof.} Let $p\in X_n$ and $q\in\P_k$, $\|q\|_2\leq 1$.
Since $\F^2(H_k)=\sum_{m=0}^\infty\oplus X_m$, we write $q$ as
$q=\sum_mr_m$, where $r_m\in X_m$. Notice that $pr_m\in X_{n+m}$ and
hence all of them are orthogonal. Moreover, if $x_1, x_2\in X_n$,
$y_1, y_2\in X_m$ and $x_1y_1=x_2y_2$, then it is necessary that
$x_1=x_2$ and $y_1=y_2$. This implies that there is no cancellation
in $pr_m$ and hence, by the previous lemma,  $\|pr_m\|_2=\|p\|_2\|r_m\|_2$. 
Therefore,
$$\|pq\|_2=\sqrt{\sum_{m=0}^\infty\|pr_m\|_2^2}=
\sqrt{\sum_{m=0}^\infty\|p\|_2^2\|r_m\|_2^2}=\|p\|_2\|q\|_2,$$
and  $\|p\|_\infty=\|p\|_2$.

The completely bounded case is very similar. A canonical element
of $M_{n_0}(X_n)$ looks like $T=\sum_{i\leq k}A_i\otimes e_{x_i}$, where
$A_i\in M_{n_0}$ and $e_{x_i}\in X_n$. A canonical element of 
$\ell_2^{n_0}(X_m)$ looks like $q=\sum_j b_j\otimes e_{y_j}$, where
$b_j\in\ell_2^{n_0}$ and $e_{y_j}\in X_m$.  Then
$Tq=\sum_i\sum_j A_ib_j\otimes e_{x_i}e_{y_j}\in \ell_2^{n_0}(X_{n+m})$
and all of those terms are orthogonal to each other. Then the
proof proceeds as those of section 2. The complementation part is
very easy: If $p\in X_n$, then $\|p\|_\infty=\|pe_0\|_2$.\qed
\medbreak

Multiplication from the left by any one of the $e_i$'s in
$\A(k)$ or $\F^\infty(k)$ is an isometry; (i.e., 
$\|p\|_\infty=\|e_ip\|_\infty$). However, multiplication from
the right does not have to be like that.

\proclaim{Proposition 17} Let $p\in\A(k-1)\subset \A(k)$. Then
$\|pe_k\|_\infty=\|pe_k\|_2=\|p\|_2$.
\endproclaim

\noindent{\bf Proof.} Let $p\in\P_{k-1}$, $p=\sum_ia_i e_{x_i}$ where
$x_i\in P_{k-1}$, and $q\in\P_k$, $q=\sum_j b_je_{y_j}$ where
$y_j\in P_k$. Then
$$pe_kq=\sum_i\sum_j a_ib_j e_{x_i}e_ke_{y_j}.$$
Since $e_{x_i}e_ke_{y_j}=e_{x_{i'}}e_ke_{y_{j'}}$ iff $x_i=x_{i'}$
and $y_j=y_{j'}$, then Lemma 15 applies and we have that 
$\|pe_kq\|_2=\|p\|_2\|e_kq\|_2=\|p\|_2\|q\|_2$.\qed
\medbreak

We conclude with the following two applications of the previous
propositions. 

1. Let $p(e_1, e_2, \cdots, e_k)\in\P$ be a non-commutative  
homogeneous polynomial of
degree $n$; i.e., $p(\lambda e_1, \cdots, \lambda e_k)=
\lambda^n p(e_1,\cdots,e_k)$, (or $p\in X_n$ 
with the notation of Proposition 16).
If $\|\sum_{i\leq k} T_iT_i^*\|\leq1$, then
$$\|p(T_1,T_2,\cdots, T_k)\|\leq\|p\|_2.$$

\smallbreak
2. Let $T_1, T_2\in B(\ell_2)$ be such 
that $\|[T_1\  T_2]\|\leq 1$ (i.e., $\|T_1 T_1^*+T_2 T_2^*\|\leq 1$)
and let $p(t)$
be a polynomial in one variable. The classical Von Neumann's
inequality states that $\|p(T_1)\|\leq \|p \|_\infty$. Therefore, 
using the Banach algebra properties of $B(\ell_2)$ we get that
$$\|p(T_1)T_2\|\leq\sup_{t\in{\Bbb T}}|p(t)|,$$
but if we apply Proposition 17 to Popescu's inequality we get
$$\|p(T_1)T_2\|\leq\sqrt{\int_{\Bbb T}|p(t)|^2dm(t)}.$$
\medbreak
\noindent{\bf Acknowledgment.} The author thanks Gelu Popescu
for useful discussions.
\medbreak

\enddocument